\documentstyle[12pt]{amsart}
\newtheorem{theorem}{Theorem}[section]
\newtheorem{corollary}[theorem]{Corollary}
\newtheorem{remark}[theorem]{Remark}
\newtheorem{lemma}[theorem]{Lemma}
\newtheorem{proposition}[theorem]{Proposition}
\newtheorem{definition}{Definition}[section]

\numberwithin{equation}{section}

\begin{document}

\title[Laplace's equations on n-dimensional singular manifolds]
{Laplace's equations on n-dimensional singular manifolds}

\author{Fangshu Wan}

\address{Fangshu Wan, School of Mathematical Sciences\\
University of Science and Technology of China\\ Hefei 230026, P. R. China.}
\email{wfangshu@@mail.ustc.edu.cn}

\keywords{Laplace's equations, Sobolev type embeddings, Singular manifolds, the Conical metric, the Poincar\'{e} like metric}

\thanks {The work is supported by NSFC, No. 11721101.}

\subjclass{Primary 58J05; Secondary 35J15}

\begin{abstract}
We show that the Sobolev embedding is compact on punctured manifolds with conical singularities.
On the other hand, we find the Sobolev inequality does not hold on punctured manifolds with the Poincar\'{e} like metric, on
which one has the Poincar\'{e} inequality. Applying the results to the Laplace's equation on the singular manifolds, we obtain
the existence of the solution in both cases. In the conical singularity case, we prove further that the solution can be extended to singular points
and it is H\"{o}lder continuous. However, the solution can not be continuously extended to singular points in the Poincar\'{e}
like metric case. Moreover, on singular manifolds with conical singularities, we obtain the existence and regularity of nontrivial nonnegative solutions for the semilinear elliptic equation with subcritical exponents.

\end{abstract}

\maketitle

\section{Introduction}
We study the existence and regularity of solutions for the Laplace's equation on Riemannian manifolds with conical singularities and with the Poincar\'{e} like metric, which will be defined as follows.

We assume $M^0$ is an open manifold of dimension $n$ $(n \geq 2)$ which comes from a compact manifold $M$ by
removing some finite points, that is, $M^0 = M \setminus \{p_1,p_2,\ldots,p_k\}.$ At each puncture $p_i$ $(i = 1,2,\ldots,k),$ we choose a neighborhood $U_i$ of $p_i$ such that $U_i \bigcap U_j = \emptyset$ $ (i \neq j)$ and $\varphi _i: U_i \rightarrow B_1(0)$
is a diffeomorphism with $\varphi_i (p_i) = 0$ where $B_1(0)$ is the unit ball. We introduce the Riemannian metric
in $U_i^* = U_i \setminus \{p_i\},$ in a local coordinate centered at $p_i,$
\begin{equation}
\label{2.1}
ds^{2}=|x|^{2\beta_{i}} ds_{0}^{2}, \; \; \; \beta_{i} >-1,
\end{equation}
where $x\in B^*_1(0) = B_1(0) \setminus \{0\}$ and $ds_0^2$ is the Euclidean metric. Then we extend this metric to $M^0$. Obviously, the manifold is noncomplete with a finite volume.

\begin{definition}
A metric on $M^0$ is called a conical metric $g$ if it is quasi-isometric to the metric we constructed above.
\end{definition}

Moreover, we can also introduce the Poincar\'{e} like metric
in $U_i^* = U_i \setminus \{p_i\},$ in the local coordinates,
\begin{equation}
ds^2 = \frac{ds_0^2}{|x|^2(1-\log|x|)^2}
\end{equation}
where $x\in B^*_1(0) = B_1(0) \setminus \{0\}$. Then we extend this metric to $M^0$. Clearly $M^0$ is a complete manifold with a finite volume.

\begin{definition}
A complete metric on $M^0$ is called a Poincar\'{e} like metric $\omega$ if it is quasi-isometric to the metric we constructed above.
\end{definition}

Let $\nabla$ and $\Delta$ be the gradient and the Laplace-Beltrami operator with the Euclidean metric.  Associated to the conical metric $g$, one can define gradient $\nabla_g$ and $\Delta_g$ operator
in the usual way. One can also define the Hilbert space $H_g^{1}(M^0)$ which is the completion of $C^\infty (M)$ with norm
$\|\nabla_g u\|_{L_g^2}+\|u\|_{L_g^2}$, where $\|u\|_{L_g^p}=(\int_{M}|u|^{p}dV_g)^{1/p}$ is the
$L^{p}$-norm with the conical metric $g$. Similarly, we can define gradient $ \nabla_\omega$, Laplace-Beltrami operator $ \Delta_\omega$, the Banach space $L^p_{\omega}(M^0)$ and the Hilbert space $H^1_\omega(M^0)$ with the Poincar\'{e} like metric.
Let $\alpha:=\min _{i}\{\beta_{i}\}$. It follows from \eqref{2.1} that
\begin{equation}
\label{2.2}
\alpha >-1.
\end{equation}

Conical singular Riemann surfaces have been extensively studied by many authors (for example, \cite{BD, BDM, CP, FL, GJ, JWZ, LZ, MP}). Troyanov \cite{MT} studied the prescribing curvature problem on the surfaces with conical singularities. Let $M$ be a compact Riemann surface with conical singularities of divisor $\beta=\sum_{i}^{k}\beta_{i}p_{i}$. He proved the Sobolev inequality, a compactness result and the following Trudinger inequality:
$$ \int_M e^{bu^2} dV_g \leq c_b, \;\;\; \forall b < b_0,$$
for all $u \in H_g^1(M^0)$ satisfying $\int_M |\nabla_g u|^2 \; dV_g \leq 1$, $\int_M u \; dV_g =0$, where $c_b$ is a constant related to $b$ and $b_0= 4 \pi \min_i \{1,1+ \beta_i \}.$ Then Chen \cite{CW} further proved a Trudinger inequality with the best constant on such surfaces. For recent development of this topic, we refer the reader to \cite{LY, Y, YZ} and the references therein.

In this paper, on n-dimensional manifolds with conical singularities we establish a Sobolev type embedding
\begin{equation}
H^1_g(M^0)\hookrightarrow L^p_g(M^0)
\end{equation}
for any $1 \leq p \leq \frac{2n}{n-2}$ with $n>2$ and show that the embedding is compact for any $1 \leq p < \frac{2n}{n-2}$. Note that the results for $n=2$ have been obtained by Troyanov in \cite{MT}. Similar Sobolev type embeddings can be found in \cite{GGW, GMWG}.

On a compact Riemann surface with a finite set of punctures, Li and Wang \cite{LW} introduced a complete Poincar\'{e} like metric near the punctures and obtain smooth solutions of the Laplace's equation with the singular metric in Riemann surfaces excepting finite punctures. Dey \cite{DR} studied the prescribed negative Gaussian curvature problem on punctured Riemann surfaces with complete metrics. We are interested in the basic analysis properties of the singular manifolds with the Poincar\'{e} like metric, on which we find that the Sobolev's inequality does not hold (even $n = 2$) by constructing a counter example. However, on the singular manifolds we can prove that the Poincar\'{e} inequality holds.

Using the results above, we study the Laplace's equation. We say that $u$ is a weak solution of the equation with the conical metric $$- \Delta_g u = f, \;\;\; f \in L_g^2 (M^0), \eqno(P)$$
if $u \in H^1_g(M^0)$ and
\begin{equation}
\label{3.3}
\int _{M^0} \nabla_g u \nabla_g \phi \; d V_g = \int_{M^0} f \phi \;d V_g, \;\;\;\; \forall \phi \in H^1_g(M^0).
\end{equation}
Similarly, we say that $u$ is a weak solution of the equation with the Poincar\'{e} like metric
$$ -\Delta_\omega u = f, \;\;\; f \in L_\omega^2 (M^0), \eqno(Q)$$
if $u \in H^1_\omega(M^0)$ and
\begin{equation}
\label{3.3-100}
\int _{M^0} \nabla_\omega u \nabla_\omega \phi \; d V_\omega = \int_{M^0} f \phi \;d V_\omega, \;\;\;\; \forall \phi \in H^1_\omega(M^0).
\end{equation}

In both cases, we prove the existence and regularities of solutions for the Laplace's equation $(P)$ and $(Q)$. In conical singularities cases, we proved it by the compact embedding given above and the standard variational method. In the Poincar\'{e} like metric cases, we show it by the Poincar\'{e} inequality and the spectrum theory. It is interesting that we find the solution in the first cases can be extended to the singularities and it is H\"{o}lder continuous at every singularity $p_i$ $(i=1,2,\dots,k)$ which can be proved by the Moser iteration, and conversely the solution in the second cases is unbounded near the punctures even if $f$ is smooth on $M^0$ which will be shown by constructing an example.

In this paper, we also consider the semilinear elliptic equation with subcritical exponents
\begin{equation}
\label{1.1}
-\Delta_g u + hu =u^p, \;\; \mbox{$ 1 < p < \frac{n+2}{n-2}$}
\end{equation}
where $h \geq m >0$ with $h \in L_g^q(M^0)$ for some $q> \frac{n}{2}$ on singular manifolds with conical singularities. We obtain the existence of nontrivial nonnegative solutions of the equation \eqref{1.1} and show that the solutions are H\"{o}lder continuous at every singularity $p_i$ $(i=1,2,\dots,k).$

The main results of this paper are the following proposition and theorems. Assume $d(x,y)$ is the Euclidean distance from $x$ to $y$ and $B_r (p_i)$ is the Euclidean ball.
\begin{proposition}
\label{p2.2}
The embedding on the singular manifold $M$ with conical singularities $$H_g^1(M^0)\hookrightarrow L_g^p(M^0)$$ is compact for all $ p\in [1,\frac{2n}{n-2}).$
\end{proposition}

\begin{theorem}
\label{c3.7}
Let $u \in H_g^1(M^0)$ be a solution of equation $(P)$ with the conical metric and $f \in L_g^q(M^0)$ for some $q >\frac{n}{2}$. Then $u \in C^s (B_r (p_i))$ for some $s \in (0,1) $ depending only on $n,q,\beta_i.$ Moreover, for any $B_r (p_i) \subset M $
$$ |u(x) - u(y)| \leq C \Big(\frac{d(x,y)}{r}\Big)^s \left\{r^{-\frac{n}{2}}||u||_{L_g^2(B_r(p_i))} + r^{2-\frac{n}{q}} ||f||_{L_g^q(B_r(p_i))}\right\}$$
for any $x,y \in B_{\frac{r}{2}}(p_i),$ where $C$ is a positive constant depending only on $n,q,\beta_i.$
\end{theorem}

\begin{theorem}
\label{t5.2}
Let $u \in H_g^1(M^0)$ be a solution of equation \eqref{1.1}
with the conical metric where $h \geq m >0$ and $h \in L_g^q(M^0)$ for some $q >\frac{n}{2}$. Then $u \in C^s (B_r (p_i))$ for some $s \in (0,1) $ depending only on $n,q,\beta_i, \|h\|_{L_g^q(M^0)}$ and $\|u\|_{H_g^1(M^0)}.$ Moreover, for any $B_r (p_i)\subset M $
$$ |u(x) - u(y)| \leq C \Big(\frac{d(x,y)}{r}\Big)^s  r^{-\frac{n}{2}}||u||_{L_g^2(B_r(p_i))} $$
for any $x,y \in B_{\frac{r}{2}}(p_i),$ where $C$ is a positive constant depending only on $n,q,\beta_i, \|h\|_{L_g^q(M^0)}$ and $\|u\|_{H_g^1(M^0)}.$
\end{theorem}

\begin{theorem}
\label{t4.2}
The equation $(Q)$ with the Poincar\'{e} like metric has a smooth solution $u \in H^1_\omega(M^0)$ for $f \in C^\infty (M^0) \bigcap L_\omega^2 (M^0)$
if and only if $\int_{M^0} f \; dV_{\omega} = 0.$
\end{theorem}

The paper is organized as follows. In section 2, we establish the Sobolev type embedding and prove that the embedding is compact on conical manifolds. In section 3, we obtain the existence of  solutions for the equation $(P)$ with the conical metric and give the proof of Theorem \ref{c3.7}. In section 4, we obtain the existence of nontrivial nonnegative solutions of the equation \eqref{1.1} and give the proof of Theorem \ref{t5.2}. In the final section, we give the proof of Theorem \ref{t4.2} and show that the solution of $(Q)$ can not be continuously extended to the singular points in the Poincar\'{e} like metric cases.

In this paper, without other special notices, we denote $C$ to be a positive constant which may change from one line to another line.

\section{The compactness of Sobolev's embedding on a manifold with conical singularities}

In this section, we prove Sobolev's inequality and a compactness result on the singular manifold $M^0$ with conical singularities.

\begin{proposition}
\label{p2.1}
There exists a constant $C$ such that for all $u\in H_g^{1}(M^0)$ and all
$p\in[1,\frac{2n}{n-2}]$ $($if $n>2$$)$, we have
\begin{equation}
\label{2.3}
\|u\|_{L_g^{p}(M^0)}\leq C \;\|u\|_{H_g^{1}(M^0)},
\end{equation}
where $C$ depends on $n, \beta_i$ $(i=1,2,\cdots,k).$
\end{proposition}

{\bf Proof.} Let $\{\Omega _{l}\}_{l=1}^N$ be a finite covering of $M$ such that $\Omega_{i}$ $(i=1,2,\ldots,k)$ contains the singularity $p_i$ respectively
and $\Omega_{j}$ $(j=k+1,k+2,\ldots,N)$ does not contain any singularity $p_i$ $(i=1,2,\ldots,k)$.
Assume that $(\Omega _l,\varphi _l)$ are the corresponding charts. Consider
$\{\alpha_l\}$ a $C^\infty$ partition of unity subordinate to the covering
$\{\Omega _{l}\}$. We need prove there exists a constant $C_l$ such that every
$C^\infty$ function $u$ on $M$ satisfies
\begin{equation}
\label{2.4}
\|\alpha _lu\|_{L_g^{p}(\Omega_l)}\leq C_l\; \|\alpha _lu\|_{H_g^{1}(\Omega _l)}.
\end{equation}
It suffices to show that \eqref{2.4} holds when $i=1,2,\ldots,k$. Note that on each
$\Omega_{j}$ $(j=k+1,k+2,\ldots,N)$ there is a smooth metric.

In fact, since $|\nabla_g (\alpha _l u)|\leq |\nabla_g u|+|u||\nabla_g \alpha_l|$,
it is easily seen that
\begin{eqnarray}
||u||_{L_g^p(M)}&\leq & \sum_{l=1}^N ||\alpha_l u||_{L_g^p(\Omega_l)} \nonumber\\
&\leq & \sum_{l=1}^N C_l\; ||\alpha_l u||_{H_g^1(\Omega_l)}\nonumber\\
&\leq & \sum_{l=1}^N C_l\;(||u||_{L_g^2}+\sup_{1\leq l\leq N}|\nabla_g \alpha_l|||u||_{L_g^2}+||\nabla_g u||_{L_g^2})\nonumber\\
&\leq &\sup_{1\leq l\leq N}C_l\;N\; \Big[(1+\sup_{1\leq l\leq N}|\nabla_g \alpha _l|)||u||_{L_g^2}+||\nabla_g u||_{L_g^2}\Big] \nonumber\\
&\leq &C\;||u||_{H_g^1(M)} \nonumber
\end{eqnarray}
where $C$ is dependent of  $M$, $n$ and $\beta_i$ $(i=1,2,\cdots,k).$ Therefore \eqref{2.3} holds for all $u \in H_g^ 1(M^0)$ by density.

Let $U_i=\varphi _i(\Omega_i)$ be a bounded domain in $\mathbb{R}$$^{n}$ with $\varphi_i(p_i)=0$. By $(1.8)$ in \cite{GGW} (also see \cite{CKN, CWQ, LCW}), we have
\begin{equation}
\label{2.6.2}
\Big( \int_{U_i} |x|^{l} |\alpha_i u|^p dx \Big)^{1/p} \leq C \;\Big(\int_{U_i} |x|^{\theta} |\nabla (\alpha_i u)|^2 dx \Big)^{1/2},\;\; \forall 1 \leq p \leq \frac{2n}{n-2},
\end{equation}
where $C$ depends on $U_i$, $\theta$ and $l,$ on the condition that
\begin{equation}
\label{2.6.3}
n+\theta >2, \;\;\; \tau := l-\theta>-2, \;\;\; \theta \geq\frac{n-2}{2}\tau.
\end{equation}
We choose
\begin{equation}
\label{2.5}
\theta=(n-2)\beta_i, \; \; \; l=n\beta_i,
\end{equation}
then it is clear that \eqref{2.6.3} is equivalent to $$\beta_i >-1.$$
We therefore have the weighted Sobolev's inequality for any $1 \leq p \leq \frac{2n}{n-2},$
\begin{equation}
\label{2.6}
\Big( \int_{U_i} |x|^{n\beta_i} |\alpha_i u|^p dx \Big)^{1/p} \leq C_i \;\Big(\int_{U_i} |x|^{(n-2)\beta_i} |\nabla (\alpha_i u)|^2 dx \Big)^{1/2},
\end{equation}
where $C_i$ depends on $n,$ $\Omega_i$ and $\beta_i,$ if $\beta_i > -1.$
Obviously, we see by \eqref{2.6}
\begin{eqnarray}
\label{2.7}
\Big(\int_{\Omega_i} |\alpha_i u|^p \; dV_g \Big)^{1/p}&\leq & C\Big(\int_{U_i} |\alpha_i u|^p |x|^{n\beta_i}dx \Big)^{1/p} \nonumber \\
&\leq& C_i \;\Big(\int_{U_i} |x|^{(n-2)\beta_i} |\nabla (\alpha_i u)|^2 dx \Big)^{1/2} \nonumber \\
&\leq & C_i\; \Big(\int_{\Omega_i} |\nabla_g (\alpha_i u)|^2 dV_g \Big)^{1/2} \nonumber \\
&\leq & C_i\;||\alpha_iu||_{H_g^1(\Omega_i)}, \;\;\; \forall i=1,2,\ldots,k.
\end{eqnarray}
Thus we obtain inequality \eqref{2.4}.  \qed

The above Proposition \ref{p2.1} implies, that for any $p\in [1,\frac{2n}{n-2}]$ the imbedding $H^1_g(M^0) \hookrightarrow L^p_g(M^0)$
is continuous, if $n>2$. Moreover, we will show that the embedding is compact for $p\in [1,\frac{2n}{n-2})$ in the following.

{\bf Proof of Proposition 1.1.} Let $(\Omega_l,\varphi_l)$ $(l=1,2,\ldots, N)$ be a finite atlas of $M$
similar as those in the proof of Proposition \ref{p2.1}, each $\Omega_l$ being homeomorphic
to a ball $B$ of $\mathbb{R}$$^n$ $(n\geq2)$. Consider $\{\alpha_l\}$ a $C^\infty$ partition of unity subordinate to the covering $\{\Omega_l\}$.

Suppose that $\{u_m\}_{m=1}^{\infty}$ is a bounded sequence in $H_g^1(M^0)$. Then
$\{\alpha_l u_m\}_{m=1}^{\infty}$ is also a bounded sequence in $H_g^1(\Omega_l)$.
If $l\neq1,2,\ldots,k$, since the metric on $\Omega_l$ is smooth, then $\{\alpha_l u_m\}_{m=1}^{\infty}$
is precompact in $L_g^p(\Omega_l)$ by the standard Kondrakov's Theorem for compact Riemannian manifolds.

It suffices to show that $\{\alpha_l u_m\}_{m=1}^{\infty}$ is precompact in $L_g^p(\Omega_l)$ when $l=1,2,\ldots,k$.
Consider the functions defined on $B$, $l=1,2,\ldots,k$ being given:
$$h_m(x)=(\alpha_l u_m)\circ\varphi_l^{-1}(x).$$
By Proposition 1.1 in \cite{GGW}, we know that if $\beta_l >-1$ the embedding
\begin{equation}
\label{2.8}
H_0^{1,(n-2)\beta_l}(B) \hookrightarrow L_{n\beta_l}^p(B)
\end{equation}
is compact for any $p \in [1,\frac{2n}{n-2})$. Recall that, in \cite{GGW} $H_0^{1,(n-2)\beta_l}(B)$
is the completion of $C_0^{\infty}(B)$ under the norm induced by the inner product
$$(u,v) = \int_{B} |x|^{(n-2)\beta_l} \nabla u \cdot \nabla v \; dx,$$
and $L_{n\beta_l}^p(B)$ $(p \geq 1)$ is the space of functions $\varphi$ such that
$$|x|^{\frac{n\beta_l}{p}} |\varphi| \in L^p(B)$$
with the norm
$$||\varphi||_{n\beta_l, p} = \Big(\int _B |x|^{n\beta_l} |\varphi|^p \; dx \Big)^{\frac{1}{p}}.$$
It is easily seen that $\{h_m(x)\}_{m=1}^\infty$ is bounded in $H_0^{1,(n-2)\beta_l}(B).$
Hence $\{h_m(x)\}_{m=1}^\infty$ is precompact in $L_{n\beta_l}^p(B)$ for $p\in [1,\frac{2n}{n-2})$. Repeating this
argument successively for $l=1,2,\ldots,N$, we may extract a subsequence
$\{ \tilde {u}_m \}$ of the sequence $\{ u_m \}$ by a standard argument,
such that $\{\alpha_l \tilde{u}_m \circ \varphi_l^{-1}\}_{m=1}^\infty$ converges in
$L_{n \beta_l}^p (B)$ for each $l$ with $p\in [1,\frac{2n}{n-2})$. Note that $\beta_l = 0$ for $l=k+1,k+2,\cdots, N.$ Assume $\{ \alpha_l \tilde{u}_m \circ \varphi_l^{-1}\} $
converges to $g_l$ as $m \rightarrow \infty $, i.e.,
\begin{equation}
\label{2.9}
\int _{\Omega_l}|\alpha _l \tilde{u}_m - g_l \circ \varphi _l|^p \; dV_g \rightarrow 0 \;\;\; as \;\;\; m \rightarrow \infty \;\;\; \mbox{for all} \;\;l,
\end{equation}
and $g_l \circ \varphi_l \in L_g^p(\Omega_l).$

We extend functions $g_l \circ \varphi _l$ to be zero in $M \setminus \Omega_l$ and set $\tilde u = \sum _{l=1}^N g_l \circ \varphi _l$.
We easily see from \eqref{2.9} that $\tilde u \in L_g^p(M^0)$ and
\begin{eqnarray}
\label{2.10}
\Big(\int _M |\tilde u_m - \tilde u|^p \; dV_g \Big)^{1/p} &=& \Big(\int_M |\sum_{l=1}^N (\alpha_l \tilde u_m - g_l \circ \varphi_l)|^p \; dV_g \Big)^{1/p} \nonumber \\
&\leq &  \sum _{l=1}^N \Big( \int_{\Omega_{l}}|\alpha_l \tilde u_m - g_l \circ \varphi _l|^p \; dV_g \Big)^{1/p}.
\end{eqnarray}
This also implies that $\tilde u _m \rightarrow \tilde u$ in $L^p_g(M)$ as $m \rightarrow \infty $
and the proof of this proposition is complete. \qed

As a direct consequence, we have the Poincar\'{e} inequality similar as the Proposition 5 in $\cite{MT}$.

\begin{corollary}
Let $\psi \in L^2_g(M^0)$ be a function such that $\int_M \psi \;dV_g \neq 0$,
then there exists a constant $C$ such that $||u||_{L_g^2(M^0)} \leq C ||\nabla_g u||_{L_g^2(M^0)}$,
for all $u \in H_g^1(M^0)$ with $\int_M u \psi \;dV_g = 0$.
\end{corollary}

\section{Existence and Regularity of solutions for Laplace's equations on conical manifolds}
In this section, we study the equation
\begin{equation}
\label{3.1}
-\Delta_g u = f
\end{equation}
on $M^0$ where $f \in L_g^2(M^0)$. The existence theorem and regularity result of solutions for the equation $\eqref{3.1}$
are obtained. Applying the compact embedding theorem proved in last section, the proof of the existence is almost standard.

\begin{theorem}
\label{t3.1}
There exists a weak solution $u \in H^1_g(M^0)$ of $\eqref{3.1}$ if and only if $\int _M f \; dV_g=0.$
The solution $u$ is unique up to a constant.
\end{theorem}

{\bf Proof.} (i) The proof of necessity is obvious. It suffices to show that $\int_{M^0} \Delta_g u dV_g=0.$ For some singularity point $p_i$ and any $r>0$, we choose a smooth cut-off function $\varphi$ satisfying $0 \leq \varphi \leq 1$,
$$ \varphi(x)=0 \mbox{ in } B_r(p_i), \;\;\; \varphi(x)=1 \mbox{ in } M \setminus B_{2r}(p_i), \;\;\; |\nabla_g \varphi| \leq \frac{C}{r^{1+\beta_i}}.$$
In fact, we define $\psi \in C^\infty (\mathbb{R})$ with $0 \leq \psi \leq 1$,
$$\psi(s)=0 \mbox{ in } [-1,1],\;\;\; \psi(s)=1  \mbox{ in $\mathbb{R}$} \setminus [-2,2] \;\;\; \mbox{and } |\psi'(s)| \leq C. $$
Set $\varphi(x)= \psi(\frac{d(x,p_i)}{r})$ where $d(x,p_i)$ is the Euclidean distance from $x$ to $p_i$. Then for $r \leq |x| \leq 2r$,
\begin{eqnarray}
|\nabla_g \varphi (x)|^2 &\leq& C |x|^{-2 \beta_i} |\nabla \varphi(x)|^2  \nonumber  \\
& \leq & C |x|^{-2 \beta_i} |\psi'|\frac{|\nabla d(x,p_i)|^2}{r^2}   \nonumber \\
& \leq & \frac{C}{r^{2(1+\beta_i)}}. \nonumber
\end{eqnarray}
Therefore, we have
\begin{eqnarray}
| \int_{M \setminus B_r(p_i)} \Delta_g u \varphi dV_g |&\leq & \int_{B_{2r}(p_i) \setminus B_r(p_i)} |\nabla_g u | |\nabla_g \varphi| dV_g  \nonumber   \\
&\leq & \frac{C}{r^{1+\beta_i}} \int_{B_{2r}(p_i) \setminus B_r(p_i)} |\nabla_g u | dV_g  \nonumber \\
&\leq& \frac{C}{r^{1+\beta_i}} [\mbox{Vol}(B_{2r} \setminus B_r)]^{\frac{1}{2}} || \nabla_g u||_{L_g^2(B_{2r} \setminus B_r)}  \nonumber \\
&\leq& C r^{\frac{(n-2)(1+\beta_i)}{2}} || \nabla_g u||_{L_g^2(B_{2r} \setminus B_r)} \rightarrow 0, \;\;\; \mbox{as } r \rightarrow 0.  \nonumber
\end{eqnarray}

Moreover, $\int_{M \setminus B_r(p_i)} \Delta_g u \varphi dV_g = \int_{M \setminus B_{2r}(p_i)} \Delta_g u dV_g + \int_{B_{2r} \setminus B_r(p_i)} \Delta_g u \varphi dV_g,$ where
\begin{eqnarray}
\int_{B_{2r} \setminus B_r(p_i)} \Delta_g u \varphi dV_g &\leq & \int_{B_{2r} \setminus B_r(p_i)} |f| |\varphi| dV_g  \nonumber \\
&\leq & ||f||_{L_g^2(B_{2r} \setminus B_r)} [\mbox{Vol}(B_{2r} \setminus B_r)]^{\frac{1}{2}}  \nonumber \\
&\rightarrow & 0 \mbox{ as } r\rightarrow 0. \nonumber
\end{eqnarray}

So $\int_{M \setminus B_{2r}(p_i)} \Delta_g u dV_g \rightarrow 0$ as $r \rightarrow 0$ and hence $\int_{M^0} \Delta_g u dV_g =0$.

(ii) Existence of $u.$ If $f \equiv 0$, the solutions of $\eqref{3.1}$ are $u \equiv $ constant.
Hence suppose $f \not \equiv 0 .$ Consider the functional $I(u)= \int _{M^0} |\nabla_g u |^2 \; dV_g.$
Define $\mu = \inf I(u)$ for all $u \in \mathcal B$, with
$\mathcal B $ $= \{ u\in H^1_g(M^0): \int _{M^0} u \; dV_g = 0, \int_{M^0} uf \; dV_g=1\}$.

It is clear that $0 \leq \mu < \infty $. Let $\{u_i\}_{i=1}^ \infty$ be a minimizing sequence in $\mathcal B$.
Thus the set $\{|\nabla_g u_i|\}_{i=1}^ \infty$ is bounded in $L_g^2(M^0).$ It follows by the Poincar\'{e}
inequality that $\{u_i\}_{i=1}^\infty$ is bounded in $H^1_g(M^0)$. By the compactness of the embedding with
conical metric in Proposition $\ref{p2.2}$, there exists a subsequence $\{u_k\}$ and $u_0\in H_g^1(M^0)$ such that $||u_k -u_0||_{L_g^2} \rightarrow 0$
and $I(u_0) \leq \mu$.

Hence $u_0 \in \mathcal B$ and $I(u_0) = \mu$. Since $u_0$ minimizes the variational problem,
there exists two constants $\beta$ and $\gamma$ such that for all $\phi \in H_g^1(M^0)$:
$$\int _{M^0} \nabla_g u_0 \nabla_g \phi \; dV_g = \beta \int _{M^0} f \phi \; dV_g + \gamma \int _{M^0} \phi \; dV_g.$$
Picking $\phi = 1$ yields $\gamma = 0.$ Choosing $\phi = u_0$ implies $\beta = \mu.$
Since $\int _{M^0} u_0 f \; dV_g = 1$, $u_0$ is not constant and $\mu = I(u_0) > 0$. Set $\tilde u = u_0 / \mu$.
Then $\tilde u $ satisfies equation $\eqref{3.1}$ weakly in $H_g^1(M^0)$.   \qed

Furthermore, we not only obtain removable singularity for solutions to the equation $\eqref{3.1}$,
but also prove the H\"{o}lder continuity of the solution. Firstly, we prove local boundedness of the solutions.

\begin{theorem}
\label{t3.2}
Suppose $u \in H_g^1(M^0)$ is a solution of equation $\eqref{3.1}$ and $f \in L_g^q(M^0)$ for some $q >n/2$.
Then we have for any ball $B_{2r}(p_i) \subset M$ and $ p>1$,
\begin{equation}
\label{3.4}
\sup_{B_r(p_i)} |u| \leq C \Big(r^{- \frac{n}{p}}||u||_{L_g^p (B_{2r}(p_i))} + r^{2-\frac{n}{q}} ||f||_{L_g^q (B_{2r}(p_i))} \Big)
\end{equation}
where $C = C(n, \beta_i, q, p).$
\end{theorem}

{\bf Proof.} Similarly as before, we have to consider the boundedness of the solution near the conical singularity. Let the diameter $r$ be the Euclidean distance from any point on $\partial B_r(p_i)$ to $p_i$ and assume $r = 1$ firstly.
For some $k >0$ and $ m>0$, set $ \bar u = u^{+} +k$ and
\begin{equation}
\label{3.5}
\bar u _m = \left \{ \begin{array}{ll} \mbox{$ \bar u$}  \;\; &\mbox{if $u < m$},\\
k+m \;\; &\mbox{if $u \geq m$}.
\end{array} \right.
\end{equation}
Then we have $\nabla_g \bar u_m = 0$ in $\{u < 0\}$ and $\{ u \geq m\},$ and $\bar u_m \leq \bar u$.
Consider the test function $$ \varphi= \eta^2 (\bar u_m^ \beta \bar u - k^{\beta +1}) \in H_g^1(B_{2}) \;\;\; \mbox{and}\;\;\; \varphi=0 \;\;\;\mbox{on $\partial B_2$},$$
for some $\beta \geq 0$ and some nonnegative function $\eta \in C_0^1 (B_{2})$. Note that $\bar u_m$ is bounded. We only integrate in the set $\{u>0\}$ note that $\varphi =0$ and $ \nabla_g \varphi = 0$
in $\{u \leq 0\}$. Note also that $k \leq \bar u$ and $\bar u_m ^\beta \bar u - k^{\beta +1} \leq \bar u_m^\beta \bar u$
for $k >0$. We have
\begin{eqnarray}
&\;& \int_{\{u >0\}} \nabla_g u \;\nabla_g \varphi \;dV_g \nonumber \\
&\geq & \int_{\{0<u<m\}} \beta \eta ^2\bar u_m^\beta |\nabla_g \bar u_m|^2 \;dV_g + \int_{\{u>0\}} \eta ^2\bar u_m^\beta |\nabla_g \bar u|^2 \;dV_g
 - \int _{\{u>0\}} 2 |\nabla_g \bar u| |\nabla_g \eta | \bar u_m^\beta \bar u \eta \;dV_g\nonumber \\
&\geq& \int_{\{0<u<m\}} \beta \eta ^2\bar u_m^\beta |\nabla_g \bar u_m|^2 \;dV_g +\frac{1}{2}\int_{\{u>0\}} \eta ^2\bar u_m^\beta |\nabla_g \bar u|^2 \;dV_g
- 2\int _{\{u>0\}} |\nabla_g \eta |^2 \bar u_m^\beta \bar u^2 \;dV_g. \nonumber
\end{eqnarray}
Hence, we obtain
\begin{eqnarray}
&\;&\int_{\{0<u<m\}} \beta \eta ^2\bar u_m^\beta |\nabla_g \bar u_m|^2 \;dV_g +\int_{\{u>0\}} \eta ^2\bar u_m^\beta |\nabla_g \bar u|^2 \;dV_g \nonumber \\
&\leq& C \Big (\int _{\{u>0\}} |\nabla_g \eta |^2 \bar u_m^\beta \bar u^2 \;dV_g +\int _{\{u>0\}} |f| k^{-1} \;k \;\eta^2 \bar u_m^\beta \bar u \;dV_g \Big)  \nonumber \\
&\leq& C \Big (\int _{\{u>0\}} |\nabla_g \eta |^2 \bar u_m^\beta \bar u^2 \;dV_g +\int _{\{u>0\}} |f| k^{-1} \;\eta^2 \bar u_m^\beta \bar u^2 \;dV_g \Big).  \nonumber
\end{eqnarray}
Choose $k= ||f||_{L_g^q}$ if $f$ is not identically zero. Otherwise choose an arbitrary $k >0$ and eventually let $k \rightarrow 0^+.$ Set $w=\bar u_m^{\frac {\beta}{2}} \bar u$. Note
$$|\nabla_g w|^2 \leq (1+\beta) \{\beta \bar u_m^\beta |\nabla_g \bar u_m|^2 + \bar u_m ^\beta |\nabla_g \bar u|^2\}.$$
Therefore, we have
\begin{eqnarray}
\int_{B_2} |\nabla_g w|^2 \eta^2 \;dV_g &\leq & C\Big ((1+\beta)\int _M w^2 |\nabla_g \eta|^2 \;dV_g + (1+\beta) \int _M |f| k^{-1} w^2 \eta^2 \;dV_g \Big) \nonumber \\
&\leq & C\Big ((1+\beta)\int _M w^2 |\nabla_g \eta|^2 \;dV_g + (1+\beta)||\eta w||_{L_g^{\frac {2q}{q-1}}}^2  \Big). \nonumber
\end{eqnarray}
By the Sobolev's inequality and the interpolation with
\begin{equation}
2< \frac{2q}{q-1} < \frac{2n}{n-2} \;\mbox {if}\; q> \frac{n}{2}, \nonumber
\end{equation}
we have
\begin{eqnarray}
|| \eta w||_{L_g^{\frac{2q}{q-1}}} & \leq & \varepsilon || \eta w||_{L_g^{\frac{2n}{n-2}}} + C(n,q) \varepsilon ^{-\frac{n}{2q-n}} ||\eta w ||_{L_g^2}   \nonumber \\
& \leq & \varepsilon ||\nabla_g (\eta w)||_{L_g^2} + \varepsilon ||\eta w||_{L_g^2} +  C(n,q) \varepsilon ^{-\frac{n}{2q-n}} ||\eta w ||_{L_g^2},   \nonumber
\end{eqnarray}
for any small $\varepsilon >0.$ We choose a small constant $\varepsilon$ so that
\begin{equation}
\int_{B_2} |\nabla_g (w \eta)|^2 \;dV_g \leq C(1+ \beta) \int_{B_2} w^2 |\nabla_g \eta|^2 \;dV_g + \int_{B_2} \eta^2 w^2 \;dV_g + C(1+\beta)^{\frac{2q}{2q-n}} \int_{B_2} \eta^2 w^2 \;dV_g. \nonumber
\end{equation}
The sobolev inequality implies
\begin{eqnarray}
\Big ( \int_{B_2} |\eta w|^{\frac{2n}{n-2}}\;dV_g \Big) ^{\frac{n-2}{n}} &\leq& C (1+ \beta)\int_{B_2} w^2 |\nabla_g \eta|^2 \;dV_g + C\int_{B_2} \eta^2 w^2 \;dV_g \nonumber  \\
&& \;\;\;\;+ C(1+\beta)^{\frac{2q}{2q-n}} \int_{B_2} \eta^2 w^2 \;dV_g \nonumber \\
&\leq & C(1+\beta)^{\frac{2q}{2q-n}} \int_{B_2} (| \nabla_g \eta|^2 + \eta ^2) w^2 \;dV_g. \nonumber
\end{eqnarray}
For any $1\leq r < R \leq 2,$ consider an $\eta \in C_0^1(B_R)$ with the property
\begin{equation}
\eta \equiv 1 \; \mbox{in $B_r$} \;\;\; \mbox{and} \;\;\; |\nabla_g \eta| \leq \frac{C}{R^{1+\beta_i}-r^{1+\beta_i}}. \nonumber
\end{equation}
In fact, let $\rho (x,p_i) $ be the conical metric distance from $x$ to $p_i.$ Assume $\rho_1=\rho(x,p_i)=\frac{1}{1+\beta_i} r^{1+\beta_i}$ for $x \in \partial B_r(p_i)$ and $\rho_2=\rho(y,p_i)=\frac{1}{1+\beta_i} R^{1+\beta_i}$ for $y \in \partial B_R(p_i)$. We construct a function $\psi \in C_0^\infty (\mathbb R)$ satisfying
$$ \psi (s)=1 \;\;\; \mbox{in } [-\frac{1}{1+\beta_i}r^{1+\beta_i}, \frac{1}{1+\beta_i}r^{1+\beta_i}],$$
$$\psi (s)=0 \;\;\; \mbox{in $\mathbb{R}$}\setminus [-\frac{1}{1+\beta_i}R^{1+\beta_i}, \frac{1}{1+\beta_i}R^{1+\beta_i}], $$
and $|\psi'(s)| \leq \frac{C}{R^{1+\beta_i}-r^{1+\beta_i}}$. Set $\eta(x)=\psi(\rho (x,p_i))$. Then $\eta(x)$ satisfies
$$ \eta(x)=1 \;\;\; \mbox{in } B_r(p_i),\;\;\; \eta(x)=0 \;\;\; \mbox{in } M \setminus B_R(p_i),\;\;\; |\nabla_g \eta(x)| \leq \frac{C}{R^{1+\beta_i}-r^{1+\beta_i}} . $$

Setting $\chi = \frac{n}{n-2}$ for $n>2$ and $\chi > 2$ for $n=2,$ we obtain
$$ \Big( \int _{B_r} w ^{2 \chi} \;dV_g \Big)^{\frac{1}{\chi}} \leq C \frac{(1+ \beta)^{\frac{2q}{2q-n}}}{(R^{1+\beta_i}-r^{1+\beta_i})^2} \int _{B_R} w^2 \;dV_g.$$
Recalling the definition of $w$, we have
$$ \Big( \int _{B_r} \bar u ^{2 \chi} \bar u_m^{\beta \chi} \;dV_g \Big)^{\frac{1}{\chi}} \leq C \frac{(1+ \beta)^{\frac{2q}{2q-n}}}{(R^{1+\beta_i}-r^{1+\beta_i})^2} \int _{B_R} \bar u^2 \bar u_m^\beta \;dV_g .$$
Set $\gamma =\beta +2.$ Then we obtain
$$\Big ( \int_{B_r} \bar u_m^{\gamma \chi}\;dV_g \Big)^{\frac{1}{\chi}} \leq C \frac{(\gamma - 1)^{\frac{2q}{2q-n}}}{(R^{1+\beta_i}-r^{1+\beta_i})^2} \int _{B_R} \bar u ^ \gamma \;dV_g.$$
By letting $m \rightarrow \infty$ we conclude that
$$||\bar u ||_{L_g^{\gamma \chi}(B_r)} \leq \Big( C \frac{(\gamma - 1)^{\frac{2q}{2q-n}}}{(R^{1+\beta_i}-r^{1+\beta_i})^2} \Big)^{\frac{1}{\gamma}} ||\bar u|| _{L_g^\gamma (B_R)},$$
where $C$ is a positive constant depending only $n,q,\beta_i$ and independent of $\gamma.$
Now, taking $p >1,$ we set $\gamma = \gamma _j = \chi^j p$ and $r_j= 1+2^{-j},$ $j=0,1,\ldots,$
so that for $j=1,2,\ldots,$
$$||\bar u||_{L_g^{\gamma_j} (B_{r_j})} \leq \Big( C(n,p,q,\beta_i) \Big)^{\frac{j-1}{\chi^{j-1}}} ||\bar u||_{L_g^{\gamma _{j-1}}(B_{r_{j-1}})}.$$
Hence, by an iteration we obtain
$$||\bar u||_{L_g^{\gamma_j} (B_{r_j})} \leq \Big( C(n,p,q,\beta_i) \Big)^{\sum\frac{j-1}{\chi^{j-1}}} ||\bar u||_{L_g^p(B_2)}.$$
Letting $j \rightarrow \infty,$ we get $$\sup _{B_1} \bar u \leq C ||\bar u|| _{L_g^p (B_2)},$$ or
$$\sup _{B_1} u ^+ \leq C(||u ^+||_{L_g^p(B_2)} + k).$$

Since $u$ must be a supersolution of equation $\eqref{3.1},$ we similarly have
$$\sup _{B_1}(u^-) \leq C (||u ^-||_{L_g^p(B_2)} + k),$$
where $u^- = \max\{-u,0\}.$
Therefore, we obtain
$$\sup _{B_1} |u| \leq C (||u||_{L_g^p(B_2)} + k). $$

Now, the inequality \eqref{3.4} can be proved by a dilation argument. In fact, in a coordinate neighborhood of $p_i$, we assume the conical metric $g(x)= |x|^{2\beta_i} g'(x)$ where $g'(x)$ is a smooth Riemannian metric.
By a simple coordinate transformation $y=\frac{x}{r}$ for $x \in B_{2r}(p_i),$ define $\tilde{g'}(y)= g'(ry)$ and $\tilde{g}(y)=g(ry)$ for $y \in B_2(p_i)$. Set
$\tilde{u}(y)= u(ry)$ for $y \in B_2(p_i)$ and $\tilde{\varphi}(y)= \varphi (ry)$ with $\tilde{\varphi} \in H^1_{ \tilde g}(B_2(p_i))$ and $\tilde{\varphi}|_{\partial B_2 (p_i)}=0.$ It is easy to see that $\nabla_g^x u(x) \nabla_g^x \varphi (x) = r^{-2} \nabla_{\tilde g}^y \tilde{u}(y) \nabla_{\tilde g}^y \tilde{\varphi}(y)$ and $dV_{\tilde{g}}^y = r^{-n} dV_g^x$. Since $u(x)$ satisfies
$$\int_{B_{2r}(p_i)}\nabla_g^x u(x) \nabla_g^x \varphi (x) dV_g^x=\int_{B_{2r}(p_i)} f(x) \varphi(x)dV_g^x    $$
for any $\varphi \in H^1_g(B_{2r})$ and $\varphi |_{\partial B_{2r}}=0$,
then we see that $\tilde{u}(x)$ satisfies $$\int_{B_2} \nabla_{\tilde g}^y \tilde{u}(y) \nabla_{\tilde g}^y \tilde{\varphi}(y) dV_{\tilde g}^y = \int_{B_2} \tilde{f}(y) \tilde{\varphi}(y) d V_{ \tilde g}^y $$ for any $\tilde{ \varphi} \in H_{ \tilde g}^1 (B_2(p_i))$ and $\tilde{\varphi}|_{\partial B_2 (p_i)}=0,$ where $\tilde{f}(y)=r^{2} f(ry)$ for any $y \in B_2(p_i).$ A direct calculation shows for some $q > \frac{n}{2}$,
$$\| \tilde{f} \|_{L_{ \tilde g}^q (B_2(p_i))} = r^{2-\frac{n}{q}} \|f\|_{L_g^q(B_{2r}(p_i))} \leq C.$$
We may apply what we just proved to $\tilde{u}$ in $B_1(p_i)$ i.e.
$$ \sup_{B_1(p_i)} |\tilde{u}| \leq C \Big( \| \tilde{u} \|_{L_{ \tilde g}^p(B_2(p_i))} + \| \tilde{f}\|_{L^q_{ \tilde g} (B_2(p_i)) } \Big).$$ Then we rewrite the result in terms of $u$ and obtain
$$ \sup_{B_r(p_i)} |u| \leq C \Big( r^{-\frac{n}{p}} \|u\|_{L_g^p (B_{2r}(p_i))} + r^{2-\frac{n}{q}} \|f\|_{L_g^q(B_{2r}(p_i))} \Big).$$
\qed

\begin{remark}
\label{R3.5}
More generally, suppose $c(x) \in L_g^q (M^0)$ for some $q > \frac{n}{2}$ and
$$||c||_{L_g^q} \leq \Lambda,$$
for some positive constant $\Lambda.$ Assume that $u \in H_g^1(M^0)$ is a subsolution of $$-\Delta_g u +c(x) u(x) = f(x)$$ in the following sense
$$\int _M \nabla_g u \nabla_g \varphi \;dV_g + \int_M c u \varphi \;dV_g \leq \int_M f \varphi  \;dV_g\;\; \mbox{for any} \;\; \varphi \in H_g^1(M^0) \;\mbox{with} \;\; \varphi \geq 0 \;\; \mbox{in} \;\; M.$$
Then if $f \in L_g^q(M^0),$ the inequality
\begin{eqnarray}
\sup_{B_r(p_i)} u^+ &\leq & C \Big(r^{- \frac{n}{p}}||u||_{L_g^p (B_{2r}(p_i))} + r^{2-\frac{n}{q}} ||f||_{L_g^q (B_{2r}(p_i))} \Big)  \nonumber
\end{eqnarray}
holds for any $p >0$ where $C$ is a positive constant depending only on $n, \Lambda, \beta_i, p, q.$ The proof is similar to that of Theorem \ref{t3.2}.
\end{remark}

Then the next result is referred to as the weak Harnack inequality. To this aim,
we need a important lemma.

Because the conical metric space satisfies doubling condition, i.e.
$$ \frac{\mbox{vol}(B_{2R})}{\mbox{vol}(B_R)} \leq C_{n, \beta_i}.$$
Indeed, we can assume that $B_R$ is a Euclidean ball of radius $R$ centered at $p_i$ $(i=1,2,\ldots, k).$ Then
\begin{eqnarray}
\mbox{vol}(B_R) &\leq & C \int_{B_R} |x|^{n \beta_i} dx \nonumber \\
&\leq & C \int_0^R \int_{S^{n-1}} r^{n \beta_i + n-1} d\sigma dr \nonumber \\
&\leq & C \frac{W_n}{n(\beta_i + 1)} R^{n(\beta_i +1)}.
\end{eqnarray}

It was proved that (see \cite{SE}) the John-Nirenberg lemma in \cite{JN} generalizes easily to accommodate doubling measures
by the Calderon-Zygmund decomposition.

\begin{lemma}
\label{L3.4}
Let $u_{x,r} = \frac{1}{\mbox{vol} ( B_r(x))} \int _{B_r (x)} u \;dV_g,$ where $B_r(x)$ is the Euclidean ball of radius $r$ centered at $x.$ Suppose
$u \in L_g^1(M^0)$ satisfies $$ \frac{1}{\mbox{vol} (B_r(x))} \int_{B_r} |u -u_{x,r}| \;dV_g \leq M_0, \; \;\mbox{ for any} \;\; B_r(x) \subset M,$$
where $M_0$ is a positive constant. Then there holds for any $B_r(x) \subset M $
\begin{equation}
\label{3.9}
\frac{1}{\mbox{vol} (B_r(x))} \int_{B_r(x)} e^{\frac{p_0}{M_0} |u-u_{x,r}|} \;dV_g\leq C,
\end{equation}
for some positive $p_0$ and $C$ depending only on $n.$
\end{lemma}

\begin{theorem}
\label{t3.3}
Suppose $u \in H_g^1(M^0)$ is a nonnegative solution of equation $\eqref{3.1}$ and $f \in L_g^q(M^0)$ for some $q >n/2$.
Then we have for any ball $B_{4r}(p_i) \subset M$ and $ 0 < p <\frac{n}{n-2}$,
\begin{equation}
\label{3.11}
r^{- \frac{n}{p}} ||u||_{L_g^p (B_{2r}(p_i))} \leq C \Big(\inf _{B_r(p_i)} u + r^{2-\frac{n}{q}} ||f||_{L_g^q (B_{2r}(p_i))}\Big),
\end{equation}
where $C = C(n, \beta_i, q, p).$
\end{theorem}

{\bf Proof.} Let the diameter $4r$ be the Euclidean distance from $x \in \partial B_{4r}(p_i) $ to $p_i.$ We only prove for $r = 1.$ Then the result about $r$ can be obtained by a dilation argument
similar to the proof of Theorem \ref{t3.2}.

{\it Step 1.} We prove that the result holds for some $p_0 >0.$

Set $\bar u = u+k >0$ for some $k >0$ to be determined and $v = \bar u^{-1}.$
For any $\varphi \in H_g^1(B_4)$ with $\varphi =0$ on $\partial B_4$ and $\varphi \geq 0$ in $B_4,$ consider $\bar u ^{-2} \varphi$
as the test function in the definition of weak solution. We have
$$\int_{B_4} \nabla_g u \frac{\nabla_g \varphi}{\bar{u}^2} dV_g -2 \int_{B_4} \nabla_g u \nabla_g \bar{u} \frac{\varphi}{\bar{u}^3}dV_g = \int_{B_4} f \frac{\varphi}{\bar{u}^2} dV_g.$$
Note $\nabla_g \bar{u}=\nabla_g u$ and $\nabla_g v=-\bar{u}^{-2} \nabla_g \bar{u}.$
Therefore,
$$\int_{B_4} (\nabla_g v \nabla_g \varphi + \tilde f v \varphi )\;dV_g \leq 0,$$
where $\tilde f = \frac{f}{\bar u}.$ Then $v$ is a nonnegative subsolution of the homogeneous equation $-\Delta_g v + \tilde f v =0$.
Choosing $k = ||f||_{L_g^q}$, we have $||\tilde f||_{L_g^q} \leq 1.$
Remark $\ref{3.5}$ implies that for any $p > 0$
$$\sup _{B_1} \bar u ^{-p} \leq C\int_{B_{2}} \bar u ^{- p}\;dV_g,$$
or,
\begin{eqnarray}
\inf_{B_1} \bar u &\geq& C \Big(\int_{B_{2}} \bar u ^{-p} \;dV_g\Big)^{-\frac{1}{p}} \nonumber \\
&=& C \Big(\int_{B_{2}} \bar u ^{-p} \;dV_g \int_{B_{2}} \bar u ^p \;dV_g \Big)^{-\frac{1}{p}} \Big( \int_{B_{2}} \bar u ^p \;dV_g \Big)^{\frac{1}{p}}, \nonumber
\end{eqnarray}
where $C$ is a positive constant depending on $n,\beta_i,p,q.$

The key point next is to show that there exists a $p_0 >0$ such that
$$\int_{B_{2}} \bar u ^{-p_0}\;dV_g \int_{B_{2}} \bar u ^{p_0}\;dV_g \leq C,$$
where $C$ is a positive constant depending on $n,\beta_i,q.$ We will show for any $B_{R} \subset B_4$
$$\int _{B_{R}} e^{p_0 |w|}\;dV_g \leq C,$$
where $w = \mbox{log} \; \bar u - \beta$ with $\beta = |B_{R}|^{-1} \int_{B_{R}} \mbox{log} \; \bar u \;dV_g.$

Consider $\bar u^{-1} \varphi$ as the test function with $\varphi \in L^\infty(B_4)\bigcap H_g^1(B_4) $, $\varphi =0$ on $\partial B_4$
and $\varphi \geq 0.$ By a direct calculation we have
\begin{equation}
\label{3.12}
\int_{B_4} |\nabla_g w|^2 \varphi \;dV_g = \int_{B_4} \nabla_g w \nabla_g \varphi \;dV_g- \int_{B_4} \tilde f \varphi \;dV_g.
\end{equation}
Replace $\varphi$ by $\varphi^2$ in $\eqref{3.12}.$ Then the Cauchy inequality implies
\begin{eqnarray}
\int _{B_4} |\nabla_g w|^2 \varphi^2 \;dV_g&=&\int _{B_4} 2 \nabla_g \omega \varphi \nabla_g \varphi dV_g - \int _{B_4} \tilde{f} \varphi^2 dV_g \nonumber   \\
&\leq & \frac{1}{2} \int _{B_4} |\nabla_g w|^2 \varphi^2 dV_g + 2\int _{B_4} |\nabla_g \varphi|^2 dV_g-\int _{B_4} \tilde{f} \varphi^2 dV_g  \nonumber  \\
&\leq& C \left\{\int_{B_4} |\nabla_g \varphi|^2\;dV_g + \int_{B_4} |\tilde f| \varphi^2 \;dV_g\right\},  \nonumber
\end{eqnarray}
where
\begin{eqnarray}
\int_{B_4} |\tilde f| \varphi ^2 \;dV_g &\leq& ||\tilde f||_{L_g^{\frac{n}{2}}} ||\varphi||_{L_g^{\frac{2n}{n-2}}}^2  \nonumber \\
&\leq&C(n,q) ||\varphi||_{L_g^{\frac{2n}{n-2}}}^2 \nonumber \\
&\leq&C(n,q,\beta_i) ||\nabla_g \varphi||_{L_g^2}^2 .  \nonumber
\end{eqnarray}
Therefore, we have
$$\int_{B_4}|\nabla_g w|^2 \varphi ^2 \;dV_g \leq C(n,q,\beta_i)\int_{B_4}|\nabla_g \varphi|^2 \;dV_g.$$
For any $B_{2R}(p_i) \subset B_4,$ choose $\psi \in C_0^\infty (B_{2}(p_i))$ with
\begin{equation}
\psi(\widetilde{x}) = \left \{ \begin{array}{ll} \mbox{$ 1$}  \;\; &\mbox{$\widetilde{x}\in B_1(p_i)$},\\
0 \;\; &\mbox{$ \widetilde{x} \in M \setminus B_2(p_i)$},
\end{array} \right. \;\;\; \mbox{and $|\nabla_g \psi(\widetilde{x})| \leq C.$}
\end{equation}
By scaling, we consider the cut-off function $$\varphi(x)=\psi (\frac{x_1}{R},\cdots,\frac{x_{n}}{R})= \psi(\widetilde{x})$$
with
$$\mbox{supp}\; \varphi \subset B_{2R}(p_i), \;\varphi \equiv 1\;\mbox{in} \; B_R(p_i),\;|\nabla_g \varphi|\leq \frac{C}{R^{1+\beta_i}}.$$
Indeed,
\begin{eqnarray}
|\nabla_g \varphi(x)|^2 &\leq &C |x|^{-2 \beta_i} |\nabla \varphi(x)|^2 \nonumber \\
&=& C |R|^{-2\beta_i} |\widetilde{x}|^{-2\beta_i} R^{-2}|\nabla \psi(\widetilde{x})|^2 \nonumber \\
&\leq & C R^{-2(1+\beta_i)} |\nabla_g \psi|^2.  \nonumber
\end{eqnarray}
Hence, we obtain
$$\int _{B_R(p_i)} |\nabla_g w|^2 \;dV_g \leq C\; R^{-2(\beta_i+1)}\mbox{vol$(B_{2R}(p_i))$}. $$
By the Poincar\'{e} inequality
$$\int_{B_R(p_i)} |w-w_{y,R}|^2 \;dV_g \leq CR^{2(\beta_i + 1)} \int_{B_R(p_i)} |\nabla_g w|^2 \;dV_g,$$
we have
\begin{eqnarray}
\frac{1}{\mbox{vol}(B_R(p_i))} \int_{B_R(p_i)} |w - w_{y,R}| \;dV_g &\leq& \frac{1}{\mbox{vol}(B_R(p_i))^{\frac{1}{2}}} \Big(\int_{B_R(p_i)} |w - w_{y,R}|^2\;dV_g \Big)^{\frac{1}{2}} \nonumber \\
&\leq & \frac{1}{\mbox{vol}(B_R(p_i))^{\frac{1}{2}}} \Big(C R^{2(\beta_i +1)}\int_{B_R(p_i)} |\nabla_g w |^2 \;dV_g \Big)^{\frac{1}{2}} \nonumber \\
&\leq& C \frac{\mbox{vol}(B_{2R}(p_i))^{\frac{1}{2}}}{\mbox{vol}(B_{R}(p_i))^{\frac{1}{2}}} \leq  C.
\end{eqnarray}
Then Lemma $\ref{L3.4}$ implies
$$\int _{B_R} e^{p_0 |w| }\;dV_g \leq C.$$

{\it Step 2.} The result holds for any positive $p < n/(n-2).$

We need to prove for any $1\leq r_1 <r_2 \leq 3$ and $0< p_2<p_1<n/(n-2),$
\begin{equation}
\label{3.13}
\Big(\int _{B_{r_1}} \bar u ^{p_1} \;dV_g\Big)^{\frac{1}{p_1}} \leq C \Big( \int _{B_{r_2}} \bar u^{p_2} \;dV_g\Big)^{\frac{1}{p_2}},
\end{equation}
where $C$ is a positive constant depending on $n,q,\beta_i,r_1,r_2,p_1$ and $p_2.$

Take $\varphi = \bar u ^{-\beta}\eta^2 $ for $\beta \in (0,1)$ as the test function and we have by Cauchy inequality
$$\int_{B_4} |\nabla_g \bar u|^2 \bar u^{-\beta -1} \eta^2 \;dV_g \leq C \Big( \frac{1}{\beta^2} \int_{B_4} |\nabla_g \eta|^2 \bar u^{1-\beta}\;dV_g + \frac{1}{\beta} \int_{B_4} \frac{|f|}{k}\eta^2 \bar u ^{1-\beta} \;dV_g\Big).$$
Set $\gamma =1-\beta \in (0,1)$ and $w=\bar u ^{\frac{\gamma}{2}}.$ Then, by H\"{o}lder inequality, interpolation inequality and Sobolev inequality we obtain
$$\int _{B_4} |\nabla_g (w \eta)|^2 \;dV_g \leq \frac{C}{(1-\gamma)^m} \int_{B_4} w^2 (|\nabla_g \eta|^2+\eta^2)\;dV_g ,$$
for some positive $m >0. $ Proposition $\ref{p2.1}$ and an appropriate choice of cut-off function imply,
with $\chi = n/(n-2),$ for any $1 \leq r <R \leq 3$
$$\Big( \int_{B_r} w^{2 \chi} \;dV_g \Big)^{\frac{1}{\chi}} \leq \frac{C}{(1-\gamma)^m} \cdot \frac{1}{(R^{1+\beta_i}-r^{1+\beta_i})^2} \int_{B_R} w^2\;dV_g,$$
or
$$\Big( \int_{B_r} \bar u^{\gamma \chi} \;dV_g \Big)^{\frac{1}{\gamma \chi}} \leq \Big(\frac{C}{(1-\gamma)^m} \cdot \frac{1}{(R^{1+\beta_i}-r^{1+\beta_i})^2}\Big)^{\frac{1}{\gamma}} \Big(\int_{B_R} \bar u^\gamma \;dV_g \Big)^{\frac{1}{\gamma}}.$$
This holds for any $\gamma \in (0,1).$  Now $\eqref{3.13}$ follows after finitely many iterations.  \qed
\begin{remark}
\label{R3.6}
Suppose $c(x) \in L_g^q (M^0)$ for some $q > \frac{n}{2}$ and
$$||c||_{L_g^q} \leq \Lambda,$$
for some positive constant $\Lambda.$ Assume that $u \in H_g^1(M^0)$ is a nonnegative supersolution
of $$-\Delta_g u +c(x) u(x) = f(x)$$
in the following sense
$$\int _M \nabla_g u \nabla_g \varphi \;dV_g + \int_M c u \varphi \;dV_g \geq  \int_M f \varphi \;dV_g \;\; \mbox{for any} \;\; \varphi \in H_g^1(M^0) \;\mbox{with} \;\; \varphi \geq 0 \;\; \mbox{in} \;\; M.$$
Then if $f \in L_g^q(M^0),$ we have for any ball $B_{4r}(p_i) \subset M$ and $ 0 < p <\frac{n}{n-2}$,
\begin{equation}
\label{3.16}
r^{- \frac{n}{p}} ||u||_{L_g^p (B_{2r}(p_i))} \leq C \Big(\inf _{B_r(p_i)} u + r^{2-\frac{n}{q}} ||f||_{L_g^q (B_{2r}(p_i))}\Big),
\end{equation}
where $C = C(n, \Lambda, \beta_i, q, p).$ We can use arguments similar to the proof of Theorem \ref{t3.3} to obtain \eqref{3.16}.
\end{remark}

Now the Moser's Harnack inequality is an easy consequence of above results.
\begin{corollary}
\label{c3.6}
Let $u \in H_g^1(M^0)$ be a nonnegative solution of equation $\eqref{3.1}$ and $f \in L_g^q(M^0)$ for some $q >n/2$.
Then we have for any ball $B_{4r}(p_i) \subset M,$
\begin{eqnarray}
\sup_{B_r(p_i)} u &\leq C& \Big( \inf_{B_r(p_i)} u + r^{2-\frac{n}{q}} ||f||_{L_g^q(B_{2r}(p_i))}\Big)  \nonumber
\end{eqnarray}
where $C$ is a positive constant depending on $n,\beta_i,q.$
\end{corollary}
By Remark \ref{R3.5} and Remark \ref{R3.6}, we also have the following Harnack inequality.
\begin{remark}
\label{R3.7}
If $u \in H_g^1(M^0)$ is a nonnegative solution of equation
\begin{equation}
\label{3.17}
-\Delta_g u + c(x)u = f(x) \;\; \mbox{in} \;\; M,
\end{equation}
and $c,f \in L_g^q(M^0)$ for some $q >n/2$,
then we have for any ball $B_{4r}(p_i) \subset M,$
\begin{eqnarray}
\sup_{B_r(p_i)} u &\leq& C \Big( \inf_{B_r(p_i)} u + r^{2-\frac{n}{q}} ||f||_{L_g^q(B_{2r}(p_i))}\Big)  \\
\end{eqnarray}
where $C$ is a positive constant depending on $n,\beta_i,q$ and $\|c\|_{L_g^q(M^0)}.$
\end{remark}

The H\"{o}lder continuity follows easily from Corollary $\ref{c3.6}.$

{\bf Proof of Theorem \ref{c3.7}.} The proof is standard. We prove the estimate for $r = 1.$ Set for $r_0 \in (0,1)$
$$ M(r_0) = \sup_{B_{r_0}(p_i)}u \;\;\; \mbox{and} \;\;\; m(r_0) = \inf_{B_{r_0}(p_i)} u.$$
Then $M(r_0) < \infty$ and $m(r_0) > -\infty.$ It suffices to show
$$w(r_0) \equiv M(r_0) - m(r_0) \leq C r_0^s \{||u||_{L_g^2(B_1)} + ||f||_{L_g^q(B_1)}\} \;\;\; \mbox{for any }\;\; r_0<\frac{1}{2}. $$
Set $\delta =(2-\frac{n}{q})(\beta_i +1)>0.$ Apply Corollary $\ref{c3.6}$ to $M(r_0) - u \geq 0$ in $B_{r_0}(p_i)$ to get
$$\sup _{B_{\frac{r_0}{2}}} (M(r_0) - u) \leq C \left\{\inf_{B_{\frac{r_0}{2}}} (M(r_0) - u) + r_0^\delta ||f||_{L_g^q(B_{r_0})}\right\},$$
i.e.,
\begin{equation}
\label{3.14}
M(r_0) - m(\frac{r_0}{2}) \leq C \left\{ (M(r_0) - M(\frac{r_0}{2})) + r_0^\delta ||f||_{L_g^q(B_{r_0})}\right\}.
\end{equation}
Similarly, apply Corollary $\ref{c3.6}$ to $u - m(r_0)\geq 0$ in $B_{r_0}(p_i)$ to get
\begin{equation}
\label{3.15}
M(\frac{r_0}{2}) - m(r_0) \leq C \left\{ (m(\frac{r_0}{2}) - m(r_0)) + r_0^\delta ||f||_{L_g^q(B_{r_0})}\right\}.
\end{equation}
Then by adding $\eqref{3.14}$ and $\eqref{3.15}$ together, we obtain
$$w(\frac{r_0}{2}) \leq \frac{C-1}{C+1} w(r_0)+ C r_0^\delta ||f||_{L_g^q(B_{r_0})}. $$
Obviously, by Lemma 8.23 in \cite{GT} we can show
\begin{equation}
\label{3.15.1}
w(r_0) \leq C r_0^s \Big ( \|u \|_{L_g^2 (B_1(p_i))} + \|f \|_{L^q_g(B_1(p_i))} \Big)
\end{equation}
for some $s \in (0,1)$ and any $r_0 \leq \frac{1}{2}.$

In the following, assume the conical metric $g(x)=|x|^{2 \beta_i} g'(x)$ where $g'(x)$ is a smooth Riemannian metric. We can obtain the desired result by a scaling $y= \frac{x}{r}$ with $x \in B_r (p_i).$
Define $\tilde{g}(y)= g(ry)=r^{2 \beta_i} |y|^{2 \beta_i} g'(ry)$ and
$\tilde{u} (y)= u(ry)$ for $y \in B_1 (p_i).$ It is easy to see that $\tilde{u} (y)$ satisfies
$$\int_{B_1(p_i)} \nabla_{\tilde g}^y \tilde{u} \nabla_{\tilde g}^y \tilde{\varphi} d V_{\tilde g}^y = \int_{B_1(p_i)} \tilde{f} \tilde{\varphi} dV_{\tilde g}^y$$
for any $\tilde{\varphi} \in H_{\tilde g}^1 (B_1(p_i)) $ and $\tilde{\varphi} |_{\partial B_1(p_i)}=0,$
where $\tilde{f}(y)= r^{2} f(ry)$ for $y \in B_1(p_i).$
We may apply \eqref{3.15.1} to $\tilde{u}$ in $B_1(p_i)$ and rewrite the result in terms of $u$
i.e. $$ w(\frac{r}{2}) \leq C r_0^s \Big (r^{-\frac{n}{2}} \|u\|_{L_g^2 (B_r(p_i))} + r^{2-\frac{n}{q}} \|f\|_{L_g^q (B_r(p_i))} \Big).$$
Hence, we have
$$ |u(x)-u(y)| \leq C \Big(\frac{d(x,y)}{r}\Big)^s \Big (r^{-\frac{n}{2}} \|u\|_{L_g^2 (B_r(p_i))} + r^{2-\frac{n}{q}} \|f\|_{L_g^q (B_r(p_i))} \Big)$$
for any $x,y \in B_{\frac{r}{2}} (p_i)$ and some $s \in (0,1),$ where $d(x,y)$ is the Euclidean distance from $x$ to $y.$

\qed

The H\"{o}lder continuity also follows easily from Remark $\ref{R3.7}.$
\begin{corollary}
\label{c3.8}
Let $u \in H_g^1(M^0)$ be a solution of equation $\eqref{3.17}$ with conical metric and $c,f \in L_g^q(M^0)$ for some $q >\frac{n}{2}$. Then $u \in C^s (B_r (p_i))$ for some $s \in (0,1) $ depending only on $n,q,\beta_i,$ $\|c\|_{L_g^q(M^0)}$ and $\|f\|_{L_g^q(M^0)}.$ Moreover, for any $B_r(p_i) \subset M $
$$ |u(x) - u(y)| \leq C \Big(\frac{d(x,y)}{r}\Big)^s \left\{r^{-\frac{n}{2}} ||u||_{L_g^2 (B_r(p_i))} + r^{2-\frac{n}{q}} ||f||_{L_g^q(B_r(p_i))}\right\}$$
for any $x,y \in B_{\frac{r}{2}}(p_i),$ where $C$ is a positive constant depending only on $n,q,\beta_i,$ $\|c\|_{L_g^q(M^0)}$ and $\|f\|_{L_g^q(M^0)}.$
\end{corollary}

\section{Existence and Regularity of solutions for semilinear elliptic equations with subcritical exponents on conical manifolds}
In this section, we consider the semilinear elliptic equation
\begin{equation}
\label{3.18}
-\Delta_g u + hu =u^p, \;\; \mbox{$ 1 < p < \frac{n+2}{n-2}$}
\end{equation}
on $M$ where $h \geq m >0$ and $h \in L_g^q(M^0)$ for some $q> \frac{n}{2}$. We allow the function $h$ to be unbounded from above. The existence theorem and regularity result of nontrivial nonnegative solutions for the equation $\eqref{3.18}$ are obtained.

Using the compact embeddings established in Proposition 1.1, we obtain existence of nontrivial nonnegative solutions to the equation $\eqref{3.18}$ via variational methods.

\begin{theorem}
\label{t5.1}
There exists a weak solution $u \in H_g^1(M^0)$ of $\eqref{3.18}.$
\end{theorem}

{\bf Proof.} Set $\mathcal{A}$ $ = \{u \in H_g^1 (M^0): ||u||^2_{L_g^{p+1}} = 1\}.$ Consider the minimizing problem $$J= \inf _{u \in \mathcal{A}} \int _M |\nabla_g u|^2 + hu^2\;dV_g.$$ On the one hand, we easily see from the assumption and the embedding theorem in Proposition \ref{p2.1} that
\begin{eqnarray}
\int _M hu^2\;dV_g &=& \Big(\int_M h^{q}\;dV_g \Big)^{\frac{1}{q}} \Big(\int_M u^{\frac{2q}{q-1}}\;dV_g \Big)^{1-\frac{1}{q}} \nonumber \\
&\leq & C (\int_M |\nabla_g u|^2 + |u^2| \;dV_g) \nonumber \\
&\leq& C. \nonumber
\end{eqnarray}
Note that $$\frac{2q}{q-1} < \frac{2n}{n-2},\;\; \mbox{if $q>\frac{n}{2}$}.$$
This implies $J < \infty.$ On the other hand, we have
\begin{eqnarray}
\int _M |\nabla_g u|^2 + hu^2 \;dV_g &\geq& \min\{1,m\} \int_M |\nabla_g u|^2 + u^2 \;dV_g \nonumber \\
&\geq& C \min\{1,m\} ||u||_{L_g^{p+1}}^2. \nonumber
\end{eqnarray}
So it implies that $J >0.$ Using compactness of the embedding given in Proposition 1.1 and the standard variational method, we know that $J$ attains at some $u_0 \in \mathcal{A}.$ Moreover, $u_0 \geq 0$ in $M.$ Note that
$$\int _M |\nabla_g u|^2 + hu^2 \;dV_g \geq \int _M |\nabla_g |u||^2 + hu^2 \;dV_g.$$
And $u_0 \not \equiv 0$ in $M.$ We easily see that $u_0$ satisfies the equation
$$- \Delta_g u +hu = \lambda u^p$$ with $$\lambda = \int_M |\nabla_g u_0|^2 + h u_0^2 \;dV_g.$$
Then $u= \lambda ^{\frac{1}{p-1}} u_0$ is a nontrivial nonnegative solution of $\eqref{3.18}.$ The regularity of elliptic equations implies that $u \in C^\infty (M \setminus \{p_1,p_2, \ldots ,p_k\})$ when $h \in C^\infty (M \setminus \{p_1,p_2, \ldots ,p_k\}).$ The strong maximum principle also implies that $u > 0$ in $M \setminus \{p_1,p_2, \ldots ,p_k\}.$ This completes the proof of Theorem $\ref{t5.1}.$                                        \qed
\begin{remark}
\label{R3.8}
The proof of Theorem $\ref{t5.1}$ implies that when $p=1,$ the eigenvalue problem
$$- \Delta_g u +hu = \lambda u\;\; \mbox{in} \; \; M$$
admits a positive eigenvalue $0< \lambda _1 < \infty$, the corresponding eigenfunction $\phi _1 \in C^\infty (M \setminus \{p_1,p_2, \ldots ,p_k\})$ and $\phi _1 >0$ in $M \setminus \{p_1,p_2, \ldots ,p_k\}$, if $h \in C^\infty (M \setminus \{p_1,p_2, \ldots ,p_k\}).$
\end{remark}

We will show that the solution $u$ of $\eqref{3.18}$ is bounded and H\"{o}lder continuous at every singularity $p_i$ $(i=1,2,\dots,k)$ in the following.

{\bf Proof of Theorem \ref{t5.2}.} We apply Corollary \ref{c3.8} to show the H\"{o}lder continuity of the solution at every singularity $p_i$ $(i=1,2,\dots,k)$ for the equation $\eqref{3.18}$. In our case, we choose $c(x) = h-u^{p-1} $ and $f(x) = 0$ in the equation $\eqref{3.17}.$ Indeed, it follows from $u \in H_g^1(M^0)$ and the embedding theorem in Proposition \ref{p2.1} that $||u^{p-1}||_{L_g^q} < \infty$ provided $q \leq \frac{2n}{(n-2)(p-1)}.$ Note that $\mbox{vol}(M^0)< \infty$. Since
$$\frac{2n}{(n-2)(p-1)} > \frac{n}{2} \;\;\; \mbox{for $1< p < \frac{n+2}{n-2}$},$$
then we can choose $$\frac{n}{2} < q < \frac{2n}{(n-2)(p-1)}$$ such that $||u^{p-1}||_{L_g^q} < \infty.$ This completes the proof of Theorem \ref{t5.2}.          \qed

\section{Laplace's equations on manifolds with the Poincar\'{e} like metric}
In this section we mainly prove the existence of the solution for $\Delta_\omega u = f$
on punctured Riemannian manifolds with the Poincar\'{e} like metric.

We denote the Laplace-Beltrami of $M^0$ induced by the Poincar\'{e} like metric by $\Delta_\omega $
which is the closed extension of Laplacian acting on smooth functions with compact support. Moreover, the definition
domain of $\Delta_\omega$ called $\mathcal{D}$$(\Delta_\omega)$ consists of those $u \in L_\omega^2(M^0)$ such that $\Delta_\omega u \in L_\omega^2(M^0)$ in
the sense of distribution. Then $\Delta_\omega$ is a self-adjoint operator and $\mathcal D$$(\Delta_\omega)$ $ \subset H^1_\omega(M^0)$ because
$||\nabla_\omega u||_2^2 = (u,\Delta_\omega u).$ Since $M^0$ has a finite volume, constants are eigenfunctions of $\Delta_\omega.$

In the following theorem, we prove the Poincar\'{e} inequality for the complete Riemannian manifold $M^0.$ Let $B_R^M(p)$ be the ball with the Poincar\'{e} like metric $\omega.$

\begin{theorem}
\label{t4.1}
For any $f \in H^1_\omega(M^0)$ with $\int_{M^0} f \;dV_{\omega} = 0,$ we have
$$\int_{M^0} f^2 \;dV_{\omega} \leq \lambda _1^{-1} \int_{M^0} |\nabla_\omega f|^2 \;dV_{\omega}.$$
\end{theorem}

{\bf Proof.} On the Reimannian manifold with the Poincar\'{e} like metric $\omega$, suppose $\lambda _1(U_i^*)$ is the first Dirichlet eigenvalue of $U_i^*.$
We claim that $\lambda _1(U_i^*) \geq C.$

We set $h(r) = r^2(1-\log r)^2$ where $r = d(x,p_i)$ is the Euclidean distance from $x$ to $p_i$. Then the metric of $U_i^*$ can be expressed by $ds^2 = \frac{ds_0^2}{h(r)}.$
The volume element $dV_{\omega} = \frac{1}{h(r)^{\frac{n}{2}}}dx,$ where $x = (x_1,x_2,\ldots,x_n).$ If we use the spherical coordinates
we have $ds^2 = \frac{1}{h(r)}(dr^2 + r^2 d \sigma^2)$ and $dV_{\omega} = \frac{r^{n-1}}{(h(r))^{\frac{n}{2}}} dr d\sigma.$ Clearly it suffices to show this statement using the quasi-isometric metric $ds^2 = \frac{ds_0^2}{h(r)}.$ 

For any $g \in C_0^\infty (U_i^*)$ we show that
\begin{equation}
\label{4.1}
\int_{U_i^*} |g| \;dV_{\omega} \leq \frac{1}{n-1} \int_{U_i^*} |\nabla_\omega g| \;dV_{\omega}.
\end{equation}

Since$$|\nabla_\omega g|^2 = h(r)\Big(\frac{\partial g}{\partial r}\Big)^2 + \frac{h(r)}{r^2}\Big(\frac{\partial g}{\partial \theta}\Big)^2 \geq h(r) \Big(\frac{\partial g}{\partial r}\Big)^2,$$
we have
$$ |\nabla_\omega g| \geq \sqrt{h(r)} \left|\frac{\partial g}{\partial r}\right| \geq -\sqrt{h(r)} \frac{\partial |g|}{\partial r}.$$
Hence we obtain
$$\int_0^1 |\nabla_\omega g| \frac{r^{n-1}}{(h(r))^{\frac{n}{2}}} \; dr \geq -\int_0^1 \sqrt{h(r)} \frac{\partial |g|}{\partial r} \frac{r^{n-1}}{(h(r))^{\frac{n}{2}}}\; dr. $$
Because of the support assumption of $g,$ applying the integration by parts with no boundary terms yields
\begin{eqnarray}
\label{4.2}
\int_0^1 |\nabla_\omega g| \frac{r^{n-1}}{(h(r))^{\frac{n}{2}}} \;dr &\geq& -\int_0^1 \frac{\partial |g|}{\partial r} \frac{r^{n-1}}{(h(r))^{\frac{n-1}{2}}} \; dr \nonumber \\
&=& \int_0^1 |g| \frac{\frac{\partial}{\partial r}\Big(\frac{r^{n-1}}{(h(r))^{\frac{n-1}{2}}}\Big)}{\frac{r^{n-1}}{(h(r))^{\frac{n}{2}}}} \frac{r^{n-1}}{(h(r))^{\frac{n}{2}}} \;dr.
\end{eqnarray}
A straightforward computation gives
\begin{equation}
\label{4.3}
\frac{\frac{\partial}{\partial r}\Big(\frac{r^{n-1}}{(h(r))^{\frac{n-1}{2}}}\Big)}{\frac{r^{n-1}}{(h(r))^{\frac{n}{2}}}} = n-1.
\end{equation}
Obviously $\eqref{4.1}$ follows from $\eqref{4.2}$ and $\eqref{4.3}.$

Applying the same argument to $|g|^2$ in place of $|g|$ we obtain
$$\int _{U_i^*} |g|^2 \;dV_{\omega} \leq \frac{2}{n-1} \int_{U_i^*} |g| |\nabla_\omega g| \;dV_{\omega}.$$
By the H\"{o}lder inequality we have
$$\int _{U_i^*} |g|^2 \;dV_{\omega} \leq \frac{4}{(n-1)^2} \int_{U_i^*} |\nabla_\omega g|^2 \;dV_{\omega},$$
which means $\lambda_1(U_i^*) \geq \frac{(n-1)^2}{4}.$ Then
$$\int _{\bigcup_{i=1}^k U_i^*} |g|^2 \;dV_{\omega} \leq \frac{4}{(n-1)^2} \int_{\bigcup_{i=1}^k U_i^*} |\nabla_\omega g|^2 \;dV_{\omega}.$$
We clearly have the first Dirichlet eigenvalue $\lambda_1(\bigcup_{i=1}^k U_i^*) \geq \frac{(n-1)^2}{4}$
of $\bigcup_{i=1}^k U_i^*$.

Since $M^0$ is a complete Riemannian manifold, then $$\mbox{inf}\; \sigma(-\Delta)|_{\bigcup_{i=1}^k U_i^*}=\lambda_1(\bigcup_{i=1}^k U_i^*).$$ Therfore,
\begin{equation}\label{5.1}
\sigma(-\Delta)|_{\bigcup_{i=1}^k U_i^*} \subset [\lambda_1(\bigcup_{i=1}^k U_i^*), +\infty] \subset [\frac{(n-1)^2}{4}, +\infty].
\end{equation}
By the Donnely's decomposition principle in \cite{DH} (Lemma 5.1), which means that the essential spectrum
of $\Delta_\omega$ does not depend on the changes of the metric in a compact domain of $M^0,$ we have
\begin{equation}\label{5.2}
\sigma _{ess}(-\Delta_\omega)|_M =\sigma _{ess}(-\Delta_\omega)|_{\bigcup_{i=1}^k U_i^*},
\end{equation}
where $\sigma _{ess}(-\Delta_\omega)|_{\bigcup_{i=1}^k U_i^*}$ is the essential spectrum of $-\Delta_\omega$ in $\bigcup_{i=1}^k U_i^*$. It follows from \eqref{5.1} and \eqref{5.2} that
$$\sigma _{ess}(-\Delta_\omega)|_M \subset [\frac{(n-1)^2}{4}, +\infty). $$
By the assumption $\int_{M^0} f dV_\omega =0$, we know that $f$ is not constant. Let $$H=\{u \in H^1_\omega (M^0): \int_{M^0} f dV_\omega =0\}.$$ Then $-\Delta_\omega$ in H has only discrete spectrun in $(0, \frac{(n-1)^2}{4})$. Also, the number of the discrete spectrum in $(0, \frac{(n-1)^2}{4})$ is finite.
This implies
$$\inf \Big\{\sigma (-\Delta_\omega) \setminus \{0\} \Big\} := \mu >0.$$
Therefore, we obtain
$$ \int_M f^2 dV_\omega \leq \frac{1}{\mu } \int_M |\nabla_\omega f|^2 dV_\omega. $$
 \qed

Now we consider the existence theorem for the equation
\begin{equation}
\label{4.5}
-\Delta_\omega u = f
\end{equation}
on $M^0$ where $f \in L_\omega^2 (M^0).$

{\bf Proof of Theorem \ref{t4.2}.} The necessity part follows easily from $\int_{M^0} \Delta_\omega u \;dV_{\omega} =0.$ For some $p \in M^0$ any $R >0,$
we choose a cut-off function $\psi$ satisfying
$$\mbox{supp}\; \psi \subset B_{2R}^M (p), \;\psi \equiv 1\;\mbox{in} \; B_R^M (p),\;|\nabla_\omega \psi|\leq \frac{C}{R}.$$
Since $\mbox{vol}(M) < C $, we have
\begin{eqnarray}
|\int _{B_{2R}^M} \Delta_\omega u \psi \;dV_{\omega}| &\leq& \int _{B_{2R}^M} |\nabla_\omega u| |\nabla_\omega \psi| \;dV_{\omega} \nonumber \\
& \leq & \frac{C}{R} \int_{B_{2R}^M \setminus {B_R}^M } |\nabla_\omega u| \;dV_{\omega}  \nonumber \\
&\leq & \frac{C}{R} ||\nabla_\omega u||_{L_\omega^2} \mbox{vol}(M) \nonumber \\
&\leq & \frac{C}{R} \rightarrow 0 \;\; \mbox{as} \; R \rightarrow \infty. \nonumber
\end{eqnarray}
On the other hand, $\int _{B_{2R}^M} \Delta_\omega u \psi \;dV_{\omega}=\int _{B_{R}^M} \Delta_\omega u \;dV_{\omega} + \int _{B_{2R}^M \setminus B_R^M} \Delta_\omega u \psi \;dV_{\omega}$, where
$$|\int _{B_{2R}^M \setminus B_R^M} \Delta_\omega u \psi \;dV_{\omega}| \leq ||f||_{L_\omega^2 (B_{2R}^M \setminus B_R^M)} \mbox{vol}(B_{2R}^M \setminus B_R^M) \rightarrow 0 \;\; \mbox{as} \; R \rightarrow \infty. $$
So $\int_{B^M_R(p)} \Delta_\omega u \;dV_{\omega} \rightarrow 0$ as $R \rightarrow \infty$ and hence $\int_{M^0} \Delta_\omega u \;dV_{\omega} =0$.

To prove the sufficiency part, we consider the closed subset $L_0^2(M^0) = \{ u\in L_\omega^2(M^0): \int_{M^0} u \;dV_{\omega}= 0\}$ of $L_\omega^2 (M^0).$ Since $\int_{M^0} f \;dV_{\omega} =0$,
then $\Delta_\omega $ is a closed operator from $\mathcal{D}$$(\Delta_\omega)$ $ \bigcap L_0^2 (M^0)$ to $L_0^2 (M^0).$ Theorem $\ref{t4.1}$ implies
that $0$ does not belong to the spectrum of $\Delta_\omega$ in $L_0^2 (M^0).$ So $0$ is in the resolvent set of $\Delta_\omega$ and $(-\Delta_\omega)^{-1}$ exists in $L_0^2 (M^0).$
Therefore, the equation \eqref{4.5} has a solution
in the distribution sense. The regularity result for elliptic equations and $f \in C^\infty (M^0)$ show that the solution is smooth. \qed

Is there the corresponding Sobolev's inequality on a punctured manifold with the Poincar\'{e} like metric, similar as Proposition $\ref{p2.1}$ ? The answer is negative. For any $p >2,$ we can construct a counter example $u \in H^1_\omega(M^0)$ and $u \notin L_\omega^p(M^0).$

Let $r=d(x,p_i)$ $(0 < r \leq 1)$ and $p > 2,$ set
$$u(r) = \frac{1}{\alpha-1} (1-\log r)^{1-\alpha}, \;\;\mbox{where}\;\; \frac{3-n}{2} <\alpha \leq 1-\frac{n-1}{p}. $$
We can check $u \in H^1_\omega(U_i^*)$ easily. Note that
\begin{eqnarray}
\int_{U_i^*} |\nabla_\omega u|^2 \;dV_{\omega} &\leq & C \int_{B^*_1(0)} |\nabla u|^2 \frac{1}{r^{n-2} (1-\log r)^{n-2}} \;dx \nonumber \\
&=& C W_{n-1} \int_0^1 |\nabla u|^2 r (1 - \log r)^{2-n} \;dr \nonumber \\
&=& C W_{n-1} \int_0^1 r^{-1} (1- \log r)^{2-n-2\alpha} \;dr  \nonumber \\
&=& \frac{C W_{n-1}}{n+2\alpha -3}, \;\; \mbox{if}\;\; \alpha > \frac{3-n}{2}. \nonumber
\end{eqnarray}
Similarly, we have that $$\int_{U_i^*} |u|^2 \;dV_{\omega} < \infty, \;\; \mbox{if}\;\; \alpha > \frac{3-n}{2}.$$
But,
\begin{eqnarray}
\int_{U_i^*} |u|^p \;dV_{\omega} &\geq & C W_{n-1} \int_0^1 r^{-1} (1- \log r)^{(1- \alpha)p - n} \;dr \nonumber \\
& =& \infty, \;\; \mbox{if}\;\; \alpha \leq 1 - \frac{n-1}{p}. \nonumber
\end{eqnarray}
Hence, we can see that the Sobolev's inequality with the Poincar\'{e} like metric does not hold.

Further, the solution of $\Delta_\omega u = f$ for $f \in C^\infty (M^0)$ can not be continuous at each puncture $p_i$ $(i=1,2,\dots,k).$ In fact, the boundedness of the solution at each puncture $p_i$ $(i=1,2,\dots,k)$ does not hold. For example,
we set
\begin{equation}
\label{4.4}
u=\log (1-\log r), \;\;\; 0 < r \leq 1.
\end{equation}
By a direct computation, we have
\begin{eqnarray}
\Delta_\omega u &=& r^2 (1-\log r)^2 \Delta u + (n-2)r^3 (1-\log r)^3 \frac{\partial}{\partial r} \Big(\frac{1}{r(1-\log  r)}\Big) \frac{\partial u}{\partial r} \nonumber \\
&=& (2-n) (1-\log r) - 1 - (n-2)\log r \nonumber \\
&=& 1-n .\nonumber
\end{eqnarray}
Obviously, we can also show that $u \in H^1_\omega (U_i^*).$ Indeed,
\begin{eqnarray}
\int _{U_i^*} |\nabla_\omega u|^2 \;dV_{\omega} &\leq & C \int_{B_1^*(0)} |\nabla  u|^2 r^{2-n}(1-\log  r)^{2-n} \; dx \nonumber \\
&=& C W_{n-1} \int_0^1 \frac{1}{r (1-\log r)^n} \;dr \nonumber \\
&=& \frac{C W_{n-1}}{n-1} < \infty, \nonumber
\end{eqnarray}
and
\begin{eqnarray}
\int_{U_i^*} u^2 \;dV_{\omega} &\leq & C \int _{B_1^*(0)} u^2 \frac{1}{r^n(1-\log r)^n} \; dx  \nonumber\\
&=&C W_{n-1} \int_0^1 \Big(\log  (1-\log r)\Big)^2 \frac{1}{r (1-\log r)^n}\; dr  \nonumber \\
&=&C \frac{2}{(n-1)^3} W_{n-1} < \infty. \nonumber
\end{eqnarray}
We choose a cut-off function $\psi \in C_0^\infty (U_i)$ and $\psi \equiv 1$ in a neighborhood of $p_i,$ and set $\tilde{u} = \psi u.$ Then we have $\Delta_\omega \tilde{u} \in C^\infty (M^0).$ But $\tilde{u}$ is not bounded near each puncture $p_i$ $(i=1,2,\dots,k).$

\end{document}